\documentclass[11pt]{article}
\usepackage{graphicx} 
\usepackage{amsmath, amscd, amssymb, latexsym, epsfig, color, amsthm, tikz-cd}
\usepackage[
backend=biber,
style=numeric,
giveninits=true,
date=year,
url=false,
doi=false,
isbn=false
]{biblatex}
\renewbibmacro{in:}{}

\usepackage{enumitem}
\setlist[enumerate,1]{label=(\arabic*)}

\newtheorem{theorem}{Theorem}[section]
\newtheorem{proposition}[theorem]{Proposition}
\newtheorem{question}[theorem]{Question}
\newtheorem{corollary}[theorem]{Corollary}
\newtheorem{conjecture}[theorem]{Conjecture}
\newtheorem{lemma}[theorem]{Lemma}

\newenvironment{numtheorem}[1]{\medskip\noindent {\bf Theorem~#1.}\em}{}
\theoremstyle{definition}
\newtheorem{definition}[theorem]{Definition}

\newtheorem{remark}[theorem]{Remark}

\numberwithin{equation}{theorem}

\DeclareMathOperator{\gin}{gin}

\newcommand{\Z}{{\mathbb Z}}
\newcommand{\kk}{{\mathbb R}}
\newcommand{\Soc}{\mathrm{soc}}

\title{New necessary conditions on $h$-vectors of simplicial balls}

\author{
  Luz Elena Grisales Gómez\\
  Department of Mathematics\\
  University of Washington\\
  Seattle, WA 98195-4350, USA\\
  \texttt{luzgg@uw.edu}
  \and
  Satoshi Murai\\
  Department of Mathematics, Faculty of Education\\
  Waseda University\\
  1-6-1 Nishi-Waseda, Shinjuku, Tokyo 169-8050, Japan\\
  \texttt{s-murai@waseda.edu}
  \and
  Isabella Novik\\
  Department of Mathematics\\
  University of Washington\\
  Seattle, WA 98195-4350, USA\\
  \texttt{novik@uw.edu}
}

\date{\today}

\begin{document}

\maketitle

\begin{abstract}
    Although the $h$-vectors of simplicial balls of dimension at most $4$ have been completely characterized, much less is known in dimensions $5$ and higher. We establish new necessary conditions on the $h$-vectors of simplicial balls in these dimensions. 
    Together with the necessary conditions of Billera and Lee and the sufficient conditions of Kolins, our results yield a complete characterization of the $h$-vectors of $5$- and $6$-dimensional simplicial balls. We also give a characterization of $h$-vectors of simplicial balls of any dimension with neighborly boundary.
\end{abstract}

\section{Introduction}

A central problem at the intersection of geometric and algebraic combinatorics is to characterize the face numbers of a given class of simplicial complexes. These numbers are usually recorded in a vector called the $f$-vector. However, it is often more convenient to work with the $h$-vector, which is obtained from the $f$-vector by an invertible linear transformation and arises naturally in the Hilbert series of the Stanley--Reisner ring of the complex.

The most celebrated result of this kind is the $g$-theorem, which characterizes the $h$-vectors of simplicial spheres. The sufficiency part of this theorem is due to Billera and Lee~\cite{BilleraLeeProofSufficency1981}; the necessity part was originally proved by Stanley~\cite{Stanley1980} for simplicial polytopes and has since been established for all simplicial spheres~\cite{Adiprasito2018,AdiprasitoPapadakisPetrotou2021,KaruXiao2023,PapadakisPetrotou2020}. The statement is remarkably simple: the $g$-vectors (with entries $g_0 := 1$ and $g_i := h_i - h_{i-1}$ for $1 \leq i \leq \lfloor d/2\rfloor$) of simplicial $(d-1)$-spheres are precisely the $M$-sequences with $1+\lfloor d/2\rfloor$ entries.

Perhaps surprisingly, much less is known about the $h$-vectors of simplicial balls. The main known necessary conditions on these vectors come from simplicial spheres. If $\Delta$ is a simplicial $(d-1)$-ball, repeatedly coning over $\Delta$ produces a sequence of simplicial balls $\Delta = \Delta^{[0]}, \Delta^{[1]}, \dots, \Delta^{[d+1]}$, each of whose boundaries is a simplicial sphere. Billera and Lee~\cite{BilleraLeeNumbersFacesPolytope1981} observed that the $g$-vectors of these simplicial spheres depend only on the $h$-vector of $\Delta$ and that, by the $g$-theorem, each must be an $M$-sequence. We call these necessary conditions the \emph{Billera--Lee conditions}.

Billera and Lee conjectured that these conditions are also sufficient. This conjecture was confirmed for simplicial balls of dimension at most $4$ by Lee and Schmidt~\cite{LeeSchmidt2011} and, independently, by Kolins~\cite{Kolins2011}. However, Kolins~\cite{Kolins2011} showed that the conjecture fails in dimension $5$ and higher. He also found a new set of sufficient conditions, which we call \emph{Kolins' conditions}, that guarantee that a vector $h$ is the $h$-vector of a simplicial ball.
He proposed a conjectural characterization of the $h$-vectors of simplicial $5$-balls based on these sufficient conditions;
see Conjecture \ref{conj:Kolins} below.

Motivated by Kolins' conjecture, we establish a new set of necessary conditions for an integer vector to be the $h$-vector of a simplicial ball. Our main result is the following.

\begin{numtheorem}{\ref{thm:main}}
Let $d\geq 6$ and let $\Delta$ be a simplicial $(d-1)$-ball with $h(\Delta)=(1,h_1,\dots, h_{d-1},0)$. Then either $h_2=0$ or there exists an integer $m \geq 0$ such that
\begin{enumerate}
\item $(1, (h_1-m)-1, h_2-(h_1-m), \max\{h_3-h_2, 0\})$ is an $M$-sequence, and
\item $(1, h_1-h_{d-1}, h_2-h_{d-2}, \ldots, h_{\lfloor (d-1)/2 \rfloor}-h_{d-\lfloor (d-1)/2 \rfloor})-(0,m,0,\dots,0)$ is an $M$-sequence.
\end{enumerate}
\end{numtheorem}

Combined with the Billera–Lee conditions and Kolins’ sufficient conditions, this result confirms Kolins’ conjecture and extends the characterization one dimension further, yielding a complete characterization of the $h$-vectors of simplicial $5$- and $6$-balls.

\begin{numtheorem}{\ref{thm:5-and-6-balls}}
Let $d\in\{6,7\}$, and let $h=(1,h_1,\ldots,h_{d-1},0)$ be a nonnegative integer vector. Then $h$ is the $h$-vector of a simplicial $(d-1)$-ball if and only if
either $h=(1,h_1,0,\dots,0)$ or there exists an integer $m\geq 0$ such that the following conditions hold:
\begin{enumerate}
\item $(1, (h_1-m)-1, h_2-(h_1-m), \max\{h_{3} - h_{2},0\})$ is an $M$-sequence;
\item if $d=6$, then  $(1, (h_1-m)-h_{5}, h_2-h_4)$ is an $M$-sequence; if $d=7$, $(1, (h_1-m)-h_{6}, h_2-h_{5},h_3-h_4)$ is an $M$-sequence; and
\item $h_{\lceil d/2 \rceil} \geq h_{\lceil d/2 \rceil +1} \geq \cdots \geq h_{d-1}$.
\end{enumerate}
\end{numtheorem}

Our proof of Theorem~\ref{thm:main} is algebraic. For a simplicial ball $\Delta$ of dimension $d-1$, we consider the Stanley--Reisner ideal $I_\Delta$, a set of $d$ generic linear forms $\Theta$, and the generic initial ideal $I$ of $I_\Delta+(\Theta)$. We then analyze the socle of the polynomial ring modulo $I$.

Using one of the key lemmas used in the proof of Theorem \ref{thm:main},
we also characterize the $h$-vectors of simplicial balls of any dimension whose boundary complexes are neighborly simplicial spheres:

\begin{numtheorem}{\ref{thm:neighborly}}
Let $d\geq 6$ and let $h=(1,h_1, \ldots, h_{d-1}, 0)$ be a nonnegative integer vector. Then $h$ is the $h$-vector of a simplicial $(d-1)$-ball whose boundary is a neighborly sphere if and only if $h = (1,0,\ldots, 0)$ or the following conditions are satisfied:
\begin{enumerate}
\item There exists a nonnegative integer $s$ such that
\[
(1, h_1-h_{d-1}, h_2-h_{d-2}, \ldots, h_{\lfloor\frac{d-1}{2}\rfloor}-h_{d-\lfloor\frac{d-1}{2}\rfloor}) = \left(1, s, \binom{s+1}{2}, \ldots, \binom{s+\lfloor\frac{d-1}{2}\rfloor-1}{\lfloor\frac{d-1}{2}\rfloor}\right), 
\]
\item $(1, h_1-1, h_2-h_1, \ldots, h_{\lfloor d/2 \rfloor-1} - h_{\lfloor d/2 \rfloor-2}, \max\{h_{\lfloor d/2 \rfloor} - h_{\lfloor d/2 \rfloor-1},0\})$ is an $M$-sequence, and
\item $h_{\lceil d/2 \rceil} \geq h_{\lceil d/2 \rceil +1} \geq \cdots \geq h_{d-1}$.
\end{enumerate}
\end{numtheorem}

The remainder of the paper is organized as follows. In Section~\ref{sec:prelim}, we review the necessary background on simplicial complexes, introduce tools and concepts needed for the proofs, and recall the Billera--Lee conditions and Kolins' results in detail. In Section~\ref{sec:neighborly}, we characterize the $h$-vectors of simplicial balls with neighborly boundary, proving Theorem~\ref{thm:neighborly}. Section~\ref{sec:main} is devoted to proving the new necessary conditions in Theorem~\ref{thm:main}. We also combine these conditions with Kolins' results to characterize the $h$-vectors of simplicial $5$- and $6$-balls, proving Theorem~\ref{thm:5-and-6-balls}.

\section{Preliminaries}\label{sec:prelim}

\subsection{Simplicial Complexes}

A \emph{simplicial complex} $\Delta$ with vertex set $[n]:=\{1,2,\dots,n\}$ is a collection of subsets of $[n]$ that is closed under inclusion and contains every singleton $\{i\}$ for $i \in [n]$. The elements of $\Delta$ are called \emph{faces}. A face that is maximal under inclusion is called a \emph{facet}. The \emph{dimension} of a face $F \in \Delta$ is $\dim F := |F|-1$, and the dimension of $\Delta$ is the maximum dimension of its faces. A simplicial complex is \emph{pure} if all its facets have the same dimension.

Let $\Delta$ be a simplicial complex of dimension $d-1$. The $i$-th face number of $\Delta$, denoted $f_i(\Delta)$, is the number of $i$-dimensional faces of $\Delta$. When the complex is clear from context, we write $f_i$ instead of $f_i(\Delta)$. The $f$-vector of $\Delta$ is the list of its face numbers, $f(\Delta) := (f_{-1}, f_0, \ldots, f_{d-1})$. The $h$-vector $h(\Delta) := (h_0,\ldots,h_d)$ of $\Delta$ is the image of the $f$-vector under an invertible linear transformation; its entries are given by
$$h_i = \sum_{j=0}^i (-1)^{i-j}\binom{d-j}{i-j}f_{j-1}.$$
The $g$-vector of $\Delta$ is $g(\Delta) := (g_0, g_1, \ldots, g_{\lfloor d/2\rfloor})$, where $g_0 := 1$ and $g_i := h_i-h_{i-1}$ for $1\leq i \leq \lfloor d/2\rfloor$.

The \emph{geometric realization} $|\Delta|$ of a simplicial complex $\Delta$ is obtained by identifying each face $\tau = \{i_0, i_1, \ldots, i_k\}$ with a geometric $k$-simplex and gluing these simplices together along their common subfaces in the natural way. A simplicial complex $\Delta$ is called a \emph{simplicial $d$-ball} (respectively, a \emph{simplicial $d$-sphere}) if its geometric realization $|\Delta|$ is homeomorphic to the $d$-dimensional ball $B^d$ (respectively, the $d$-dimensional sphere $S^d$). In particular, simplicial balls and simplicial spheres are pure simplicial complexes. If $\Delta$ is a simplicial $d$-ball, then $\partial\Delta$ denotes the \emph{boundary complex} of $\Delta$: the simplicial $(d-1)$-sphere generated by the $(d-1)$-dimensional faces of $\Delta$ that are contained in exactly one facet of $\Delta$.

A simplicial complex $\Delta$ with vertex set $[n]$ is \emph{$i$-neighborly} if every $i$-element subset of $[n]$ is a face of $\Delta$. Equivalently, $\Delta$ is $i$-neighborly if $f_{i-1}(\Delta)=\binom{n}{i}$. A simplicial $(d-1)$-sphere is said to be \emph{neighborly} if it is $\lfloor d/2\rfloor$-neighborly.

\subsection{$M$-sequences}

A \emph{standard graded $\kk$-algebra} is a quotient $R = \kk[x_1,\dots,x_n]/I$ of a polynomial ring by a proper homogeneous ideal $I$, graded by degree; we write $R_i$ for its degree-$i$ component.
For a standard graded $\kk$-algebra $R$, its \emph{Hilbert function} is a function $\mathrm{hilb}_R:\mathbb Z_{\geq 0} \to \mathbb Z_{\ge 0}$ defined by
\[
\mathrm{hilb}_R(i)= \dim_\kk R_i.
\]

\begin{definition}
A sequence of integers $(h_0, h_1, \dots, h_m)$ is an \emph{$M$-sequence} if there exists a standard graded $\kk$-algebra $R$ such that $\dim_{\kk} R_i = h_i$ for all $0 \le i \le m$.
\end{definition}

\noindent In other words, $M$-sequences are the initial segments of Hilbert functions of standard graded algebras.

Macaulay~\cite{Macaulay1927} gave a purely numerical characterization of $M$-sequences. To state it, let $a$ and $k$ be positive integers. There are unique integers $a_k > a_{k-1} > \dots > a_j \ge j \ge 1$ such that
\[
  a = \binom{a_k}{k} + \binom{a_{k-1}}{k-1} + \dots + \binom{a_j}{j}.
\]
This expression is called the \emph{$k$-th Macaulay expansion} of $a$. Given this expansion, define
\[
  a_{\langle k\rangle} := \binom{a_k-1}{k-1} + \binom{a_{k-1}-1}{k-2}
  + \dots + \binom{a_j-1}{j-1}
  \quad\text{and}\quad 0_{\langle k\rangle} := 0.
\]
Macaulay's theorem can now be stated as follows.

\begin{theorem}[Macaulay~\cite{Macaulay1927}]
\label{thm:Macaulay}
A sequence of integers $(h_0, h_1, \dots, h_m)$ is an $M$-sequence if and only if $h_0 = 1$, each $h_i$ is nonnegative, and
\[
  (h_{i+1})_{\langle i+1\rangle} \le h_i \qquad \text{for } 1 \le i \le m-1.
\]
\end{theorem}

\subsection{The Stanley--Reisner ring}\label{subsec:face-ring}

Let $\Delta$ be a simplicial complex of dimension $d-1$ on the vertex set $[n]$. Let $S:=\kk[x_1,\dots,x_{n'}]$ for some $n'\ge n$. For a subset $F\subseteq[n']$, write $x_F:=\prod_{i\in F}x_i$ for the corresponding squarefree monomial, with the convention $x_\emptyset=1$. The \emph{Stanley--Reisner ideal} of $\Delta$ (in $S$) is
\[
  I_\Delta:=\bigl(x_F : F\subseteq[n'],\ F\notin\Delta\bigr),
\]
the ideal of $S$ generated by the squarefree monomials corresponding to the nonfaces of $\Delta$. The \emph{Stanley--Reisner ring} of $\Delta$ is $\kk[\Delta]:=S/I_\Delta$.
The ring $\kk[\Delta]$ is independent of the choice of $n'$ up to isomorphism. Indeed, since $\{j\}\notin\Delta$ for $n<j\le n'$, we have
$I_\Delta=\bigl(x_F : F\subseteq[n],\ F\notin\Delta\bigr)+(x_{n+1},\dots,x_{n'})$, and hence
\[
  \kk[\Delta]\cong\kk[x_1,\dots,x_n]/\bigl(x_F : F\subseteq[n],\ F\notin\Delta\bigr).
\]
We allow $n'\ge n$ because, at several points in this paper, we need to work in a polynomial ring with more variables than there are vertices of $\Delta$.

The ring $\kk[\Delta]$ is graded and has Krull dimension $d$. Its Hilbert series (the generating series of its Hilbert function) is
\begin{equation}\label{eq:hilb-face-ring}
  \operatorname{Hilb}(\kk[\Delta],t)=\frac{h_0+h_1t+\cdots+h_dt^d}{(1-t)^d},
\end{equation}
where $(h_0,\dots,h_d)$ is the $h$-vector of $\Delta$; see \cite[Sec.~II.1]{Stanley1996}.

A \emph{linear system of parameters} (l.s.o.p.) for $\kk[\Delta]$ is a sequence of $d$ linear
forms $\theta_1,\dots,\theta_d\in\kk[\Delta]_1$ such that $\kk[\Delta]/(\theta_1,\dots,\theta_d)$ is a finite-dimensional $\kk$-vector space. Since $\kk$ is infinite, a generic choice of $d$ linear forms gives an l.s.o.p.\ for $\kk[\Delta]$.
 
We say that $\Delta$ is \emph{Cohen--Macaulay (over $\kk$)} if some (equivalently, every) l.s.o.p.~$\Theta$ for $\kk[\Delta]$ is a regular sequence; that is, $\theta_i$ is a non-zero-divisor on
$\kk[\Delta]/(\theta_1,\dots,\theta_{i-1})$ for each $1\le i\le d$. Reisner \cite{Reisner1976} characterized Cohen--Macaulay complexes in terms of the homology of links of faces. As a consequence, every simplicial ball and every simplicial sphere is Cohen--Macaulay.
 
The importance of Cohen--Macaulay complexes for face enumeration stems from the following
observation: if $M$ is a graded module and $\theta$ is a linear form that is a non-zero-divisor
on $M$, then $\operatorname{Hilb}(M/\theta M,t)=(1-t)\operatorname{Hilb}(M,t)$. Applying this
$d$ times and using \eqref{eq:hilb-face-ring} yields the following result.
 
\begin{theorem}[Stanley]
\label{thm:hilb-CM}
If $\Delta$ is Cohen--Macaulay over $\kk$ and $\Theta = (\theta_1, \ldots, \theta_d)$ is an l.s.o.p.\ for $\kk[\Delta]$, then
\[
  \operatorname{Hilb}(\kk[\Delta]/(\Theta),t)=h_0+h_1t+\cdots+h_dt^d.
\]
In particular, $\dim_\kk(\kk[\Delta]/(\Theta))_i=h_i(\Delta)$ for $0\le i\le d$, and $(\kk[\Delta]/(\Theta))_i=0$ for $i>d$.
\end{theorem}

Simplicial spheres satisfy a much stronger result, known as the algebraic $g$-theorem.

\begin{theorem}[Algebraic $g$-theorem \cite{Adiprasito2018,AdiprasitoPapadakisPetrotou2021, KaruXiao2023,PapadakisPetrotou2020}]
\label{thm:algebraic-g}
Let $\Delta$ be a simplicial $(d-1)$-sphere. Then, for a generic l.s.o.p.\ $\Theta$ for $\kk[\Delta]$ and a generic linear form $\ell$, the multiplication map
\[
  \cdot\ell:(\kk[\Delta]/(\Theta))_{i-1}\longrightarrow(\kk[\Delta]/(\Theta))_i
\]
is injective for $1\le i\le\lceil d/2\rceil$ and surjective for
$\lfloor d/2\rfloor+1\le i\le d$. In particular, $\dim_\kk (\kk[\Delta]/(\Theta, \ell))_i = g_i(\Delta)$ for $0\le i\le\lfloor d/2\rfloor$, and so $g(\Delta)$ is an $M$-sequence.
\end{theorem}

\subsection{Known conditions on the $h$-vectors of simplicial balls}

Next, we review what is known about the $h$-vectors of simplicial balls. In contrast to simplicial spheres, whose $h$-vectors are characterized by the $g$-theorem, no complete characterization of the $h$-vectors of simplicial balls is known in general.

By examining the relationship between the $h$-vector of a simplicial ball and that of its boundary, Billera and Lee used the $g$-theorem to derive necessary conditions for the $h$-vector of a simplicial ball. To describe these conditions, let $\Delta$ be a simplicial $(d-1)$-ball, and let $u_1,\ldots,u_{d+1}$ be new vertices. For $0\le k\le d+1$, let
\[
  \Delta^{[k]}:=\Delta*u_1*\cdots*u_k
\]
denote the simplicial $(d+k-1)$-ball obtained by coning $\Delta$ $k$ times, with $\Delta^{[0]}=\Delta$.

\begin{definition} Let $h=(h_0,\dots,h_d)$ be an integer vector, and let $0\leq k \leq d+1$. Define 
    \[G^{[k]}(h):=
        \bigl(h_i-h_{d+k-i}\bigr)_ {0\le i\le \lfloor(d+k-1)/2\rfloor},
        \] 
    where  $h_i = 0$ for $i > d$.
    \end{definition}
    
When $h$ is fixed or clear from context, write $G^{[k]}$ for $G^{[k]}(h)$. As Billera and Lee observed,
$g(\partial \Delta^{[k]})=G^{[k]}(h(\Delta))$
for all $0\leq k\leq d+1$. Since $\partial \Delta^{[k]}$ is a simplicial $(d+k-2)$-sphere, the $g$-theorem gives the following necessary conditions.

\begin{theorem}[Billera--Lee~\cite{BilleraLeeNumbersFacesPolytope1981}] \label{thm:BilleraLee}
    Let $\Delta$ be a simplicial $(d-1)$-ball with $h(\Delta)=h$. Then $G^{[k]}(h)$ is an $M$-sequence for every $0\le k\le d+1$.
\end{theorem}

Of particular importance is the following immediate consequence:

\begin{corollary}\label{cor:fromBilleraLee}
Let $d\geq 4$ and let $\Delta$ be a simplicial $(d-1)$-ball with $h(\Delta)=h$. Then
\begin{enumerate}
\item $(1, h_1-h_{d-1}, h_2-h_{d-2}, \ldots, h_{\lfloor (d-1)/2 \rfloor}-h_{d-\lfloor (d-1)/2 \rfloor})$ is an $M$-sequence, and
\item $h_{\lceil d/2 \rceil} \geq h_{\lceil d/2 \rceil +1} \geq \cdots \geq h_{d-1}$.
\end{enumerate}
\end{corollary}

We refer to the conditions in Theorem~\ref{thm:BilleraLee} as the \emph{Billera--Lee conditions}. Lee and Schmidt \cite{LeeSchmidt2011} and, independently, Kolins \cite{Kolins2011} showed that these conditions completely characterize the $h$-vectors of simplicial balls of dimensions $3$ and $4$. Kolins \cite{Kolins2011} also showed that they are not sufficient in dimension $5$ or higher and established the following sufficient conditions:

\begin{theorem}[{Kolins~\cite[Theorems 14 and 17]{Kolins2011}}] \label{thm:Kolins}
	Let $d\geq 6$ and let $h= (1, h_1, \ldots, h_{d-1}, 0)$ be a nonnegative integer vector.  If $h$ satisfies the following conditions, then it is the $h$-vector of a simplicial $(d-1)$-ball:
	\begin{enumerate}
	\item $(1, h_1-1, h_2-h_1, \ldots, h_{\lfloor d/2 \rfloor-1} - h_{\lfloor d/2 \rfloor-2}, \max\{h_{\lfloor d/2 \rfloor} - h_{\lfloor d/2 \rfloor-1},0\})$ is an $M$-sequence, 
	\item $(1, h_1-h_{d-1}, h_2-h_{d-2}, \ldots, h_{\lfloor (d-1)/2 \rfloor}-h_{d-\lfloor (d-1)/2 \rfloor})$ is an $M$-sequence, and
	\item $h_{\lceil d/2 \rceil} \geq h_{\lceil d/2 \rceil +1} \geq \cdots \geq h_{d-1}$.
	\end{enumerate} 
\end{theorem}

We call these \emph{Kolins' conditions}. If a simplicial $(d-1)$-ball $\Delta$ has an $h$-vector satisfying these conditions, repeatedly stacking $(d-1)$-simplices onto $\Delta$ produces balls with $h$-vectors $h(\Delta)+(0,m,0,\ldots,0)$ for any nonnegative integer $m$. This gives the following corollary.

\begin{corollary}
\label{cor:Kolins}
Let $d\geq 6$ and let $h=(1,h_1,\ldots,h_{d-1},0)$ be a nonnegative integer vector. If there exists an integer $m\ge 0$ such that $h-(0,m,0,\ldots,0)$ satisfies the three conditions in Theorem~\ref{thm:Kolins}, then $h$ is the $h$-vector of a simplicial $(d-1)$-ball.
\end{corollary}

Kolins conjectured that, for $d=6$, these conditions are also necessary apart from the case of stacked balls, giving a complete characterization of the $h$-vectors of simplicial $5$-balls.\footnote{We have slightly modified Kolins’ original formulation of the conjecture to include the elementary case of stacked balls.}

\begin{conjecture}[{Kolins~\cite[Conjecture 18]{Kolins2011}}]  \label{conj:Kolins}
     Let $h := (1, h_1,\ldots, h_5, 0)$ be a nonnegative integer vector. Then $h$ is the $h$-vector of a simplicial $5$-ball if and only if either $h = (1, h_1, 0, \ldots, 0)$ or there exists an integer $m \geq 0$ such that $h - (0,m,0,\ldots,0)$ satisfies the three conditions in Theorem~\ref{thm:Kolins}.
\end{conjecture}
\noindent This conjecture motivates the results in the following sections.

\subsection{Generic initial ideals and socles} \label{sec:gin}

Let $S=\kk[x_1,\dots,x_n]$, and let $\tau$ be a term order with $x_1> x_2>\cdots>x_n$. For a homogeneous ideal $I\subseteq S$, its \emph{initial ideal} $\operatorname{in}_\tau(I)$ is the ideal generated by the leading monomials $\operatorname{in}_\tau(f)$ of all nonzero $f\in I$. The group $\mathrm{GL}_n(\kk)$ acts on $S_1$ by linear changes of coordinates, and this action extends to $S$ by graded $\kk$-algebra automorphisms. For $g\in\mathrm{GL}_n(\kk)$, write $gI=\{g(f):f\in I\}$.

\begin{theorem}[Galligo~\cite{Galligo1974}; Bayer--Stillman~\cite{BayerStillman1987}]\label{thm:gin-exists}
For any ideal $I\subseteq S$, there exist a nonempty Zariski open subset $U\subseteq\mathrm{GL}_n(\kk)$ and a monomial ideal $J$ such that $\operatorname{in}_\tau(gI)=J$ for all $g\in U$.
\end{theorem}

\noindent The ideal $J$ in Theorem~\ref{thm:gin-exists} is called the \emph{generic initial ideal} of $I$ with respect to $\tau$, denoted $\operatorname{gin}_\tau(I)$. 

For every homogeneous ideal $I$, the quotients $S/I$ and $S/\operatorname{in}_\tau(I)$ have the same Hilbert function (see \cite[Cor.~6.1.5]{HerzogHibi2011}). Also, a linear change of coordinates $g$ induces a graded isomorphism $S/I\cong S/gI$. Together, these facts give the following result.

\begin{proposition}\label{prop:gin-hilbert}
For every homogeneous ideal $I \subseteq S$ and every term order $\tau$, the algebras
$S/I$ and $S/\operatorname{gin}_\tau(I)$ have the same Hilbert function.
\end{proposition}

A monomial ideal $J\subseteq S$ is said to be \emph{strongly stable} if, for every monomial $m\in J$ and every pair of indices $i<j$ with $x_j\mid m$, the monomial $x_im/x_j$ also lies in $J$. 
A fundamental feature of generic initial ideals is that they are always strongly stable.

\begin{theorem}[Galligo~\cite{Galligo1974}; Bayer--Stillman~\cite{BayerStillman1987}]\label{thm:gin-Borel-fixed} 
For every ideal $I\subset S$, the ideal $\operatorname{gin}_\tau(I)$ is strongly stable.
\end{theorem}

From now on, let $\tau$ be the reverse lexicographic (revlex) order, and write $\operatorname{gin}(I):=\operatorname{gin}_\tau(I)$. The next result says that, numerically, multiplication by a generic linear form on $S/I$ behaves like multiplication by $x_n$ on $S/\operatorname{gin}(I)$ (see \cite[Thm.~4.3.6]{HerzogHibi2011}).

\begin{proposition}\label{prop:gin-multiplication}
Let $I \subseteq S$ be a homogeneous ideal, let $\ell \in S_1$ be a generic linear form, and let $i \ge 0$. Then the maps
\[
\cdot \ell \colon (S/I)_i \to (S/I)_{i+1}
\qquad\text{and}\qquad
\cdot x_n \colon (S/\gin(I))_i \to (S/\gin(I))_{i+1}
\]
have kernels of the same dimension.
\end{proposition}

For a finitely generated graded $S$-module $M$, the \emph{socle} of $M$ is the graded
submodule
\[
  \Soc(M):=
  \{u\in M : x_iu=0 \text{ for all } i\in[n]\}.
\]
Thus $\Soc(M)$ consists of the elements of $M$ that are annihilated by every variable. In particular, $SV=V$ for every $\kk$-subspace $V\subseteq\Soc(M)$. Consequently, if $M=S/I$ for a homogeneous ideal $I\subseteq S$, then every $\kk$-subspace $V\subseteq\Soc(S/I)$ is an ideal of $S/I$.
If $I$ is a strongly stable ideal, then for every $u\in S$, $x_nu\in I$ implies $x_iu\in I$ for every $i\in[n]$. This gives the following simple description of $\Soc(S/I)$.

\begin{lemma}
\label{sstablesocle}
If $I \subset S$ is a strongly stable ideal, then
\[
\Soc(S/I)= \{u \in S/I : x_n u =0\}.
\]
\end{lemma}

\section{The $h$-vector of simplicial balls with neighborly boundary}\label{sec:neighborly}

The goal of this section is to characterize the $h$-vectors of simplicial balls whose boundaries are neighborly spheres. We do so by showing that Kolins’ sufficient conditions are also necessary for such balls.

Throughout this section, let $\Delta$ be a simplicial $(d-1)$-ball with vertex set $[n]$, $\Delta^{[1]}=\Delta*\{n+1\}$,  and set $S:=\kk[x_1,\dots,x_{n+1}]$. With this choice of $S$, the Stanley--Reisner rings of $\Delta$, $\partial\Delta$, $\partial\Delta*\{n+1\}$, and $\partial\Delta^{[1]}$ are all quotients of $S$, and hence $S$-modules. We begin with the following definition:

\begin{definition}\label{def:K}
    Let $\Delta$ be a simplicial $(d-1)$-ball. Let $\theta_1,\ldots,\theta_d,\ell\in S_1$ be generic linear forms, and set $\Theta=(\theta_1,\ldots,\theta_d)$. Define
    \[
    K_i(\Delta):=\ker\bigl((S/(I_\Delta+\Theta))_i\xrightarrow{\cdot\ell}(S/(I_\Delta+\Theta))_{i+1}\bigr)
    \]
    and
    \[
    K_i(\partial\Delta):=\ker\bigl((S/(I_{\partial\Delta}+\Theta))_i\xrightarrow{\cdot\ell}(S/(I_{\partial\Delta}+\Theta))_{i+1}\bigr).
    \]
For simplicity, we suppress the dependence on $\Theta$ and $\ell$ in the notation.
\end{definition}

Note that $\Theta$ consists of $d$ generic linear forms. Since $\kk[\Delta]$ has Krull dimension $d$, $\Theta$ is an l.s.o.p.\ for $\kk[\Delta]$. By contrast, $\kk[\partial\Delta]$ has Krull dimension $d-1$, so $\theta_1,\ldots,\theta_{d-1}$ form an l.s.o.p.\ for $\kk[\partial\Delta]$, while $\theta_d$ acts as an additional generic linear form on $\kk[\partial\Delta]/(\theta_1,\ldots,\theta_{d-1})$.

The spaces $K_1(\Delta)$ and $K_1(\partial\Delta)$ have a natural interpretation in rigidity theory. The coefficients of $\theta_1,\ldots,\theta_d$ give a generic placement of the vertices in $\kk^d$. By a result of Lee \cite{Lee1996}, $\dim_\kk K_1(\Delta)$ is the dimension of the space of nontrivial infinitesimal motions of the resulting framework of the graph of $\Delta$. In particular, $K_1(\Delta)=0$ if and only if the graph of $\Delta$ is generically $d$-rigid. When $d\geq 4$, the same statements hold for $\partial\Delta$.

With this definition in hand, we state and prove the following lemma, which shows how certain algebraic properties of $\Delta$ are controlled by $\partial\Delta$. The lemma is a key ingredient in the proof of Theorem~\ref{thm:neighborly} and in the proof of Theorem~\ref{thm:main} in the next section.

\begin{lemma}\label{lem:K_i} Let $\Delta$ be a simplicial $(d-1)$-ball. Then $\dim_\kk K_i(\Delta) \leq \dim_\kk K_i(\partial\Delta)$ for all $0\leq i < \lfloor \frac{d-1}{2}\rfloor$.
\end{lemma}
\begin{proof}
    Recall that $S=\kk[x_1,\ldots,x_{n+1}]$, where $[n]$ is the vertex set of $\Delta$. Since $\partial\Delta^{[1]}=\Delta\cup(\partial\Delta*\{n+1\})$ and $\Delta\cap(\partial\Delta*\{n+1\})=\partial\Delta$, the Mayer--Vietoris sequence gives an exact sequence of $S$-modules
    \begin{align*}
    0 \to S/I_{\partial \Delta^{[1]}} \to S/I_{\Delta} \oplus S/I_{\partial\Delta*\{n+1\}} \to S/I_{\partial \Delta} \to 0.
    \end{align*}
    
   Let $\Theta=(\theta_1,\ldots,\theta_d)$ be the ideal generated by $d$ generic linear forms. Applying the right-exact functor $-\otimes_S S/\Theta$, using $(S/I)\otimes_S S/\Theta\cong S/(I+\Theta)$, and taking degree-$i$ components yields an exact sequence
    $$
    \bigl(S/(I_{\partial \Delta^{[1]}}+\Theta)\bigr)_i \xrightarrow{\psi} \bigl(S/(I_{\Delta}+\Theta)\bigr)_i \oplus \bigl(S/(I_{\partial\Delta*\{n+1\}}+\Theta)\bigr)_i \xrightarrow{\phi} \bigl(S/(I_{\partial \Delta}+\Theta)\bigr)_i \to 0.
    $$

Since $\partial\Delta^{[1]}$, $\Delta$, and $\partial\Delta*\{n+1\}$ are Cohen--Macaulay complexes of dimension $d-1$, and coning preserves the $h$-vector, the left and middle terms of the displayed sequence have dimensions $h_i(\partial\Delta^{[1]})$ and $h_i(\Delta)+h_i(\partial\Delta)$, respectively. Moreover, $\partial\Delta$ is a simplicial $(d-2)$-sphere, so the algebraic $g$-theorem implies that the right term has dimension $g_i(\partial\Delta)$ for $0\le i\le\lfloor(d-1)/2\rfloor$. 
Since 
\[
f_i(\partial\Delta^{[1]})=f_i(\Delta)+f_i(\partial\Delta*\{n+1\})-f_i(\partial\Delta),
\]
using the relation between the $f$- and $h$-vectors, we obtain
\begin{equation}  \label{eq:dim-equality}
h_i(\partial\Delta^{[1]})=h_i(\Delta)+h_i(\partial\Delta)-g_i(\partial\Delta)  \qquad \text{for } 0\le i\le\lfloor(d-1)/2\rfloor.
\end{equation}
Since the displayed sequence is exact at the middle and right terms, \eqref{eq:dim-equality} shows that $\psi$ is injective for $0\le i\le\lfloor(d-1)/2\rfloor$. Thus, for $0\le i\le\lfloor(d-1)/2\rfloor$, we have a short exact sequence
\[
0\to \bigl(S/(I_{\partial\Delta^{[1]}}+\Theta)\bigr)_i \xrightarrow{\psi} \bigl(S/(I_{\Delta}+\Theta)\bigr)_i \oplus \bigl(S/(I_{\partial\Delta*\{n+1\}}+\Theta)\bigr)_i \xrightarrow{\phi} \bigl(S/(I_{\partial\Delta}+\Theta)\bigr)_i \to 0.
\]

 Now, let $\ell$ be an additional generic linear form. Consider the multiplication map $\cdot \ell$. Assuming that $i< \lfloor \frac{d-1}{2}\rfloor$, we obtain the following commutative diagram:
    \[\begin{tikzcd}[column sep=tiny]
	0 & {\bigl(S/(I_{\partial\Delta^{[1]}}+\Theta)\bigr)_i} & {\bigl(S/(I_{\Delta}+\Theta)\bigr)_i \oplus \bigl(S/(I_{\partial\Delta*\{n+1\}}+\Theta)\bigr)_i} & {\bigl(S/(I_{\partial\Delta}+\Theta)\bigr)_i} & 0 \\
	0 & {\bigl(S/(I_{\partial\Delta^{[1]}}+\Theta)\bigr)_{i+1}} & {\bigl(S/(I_{\Delta}+\Theta)\bigr)_{i+1} \oplus \bigl(S/(I_{\partial\Delta*\{n+1\}}+\Theta)\bigr)_{i+1}} & {\bigl(S/(I_{\partial\Delta}+\Theta)\bigr)_{i+1}} & 0
	\arrow[from=1-1, to=1-2]
	\arrow[from=1-2, to=1-3]
	\arrow["{\cdot\ell=a}", from=1-2, to=2-2]
	\arrow[from=1-3, to=1-4]
	\arrow["{\cdot\ell=b_1\oplus b_2 }", from=1-3, to=2-3]
	\arrow[from=1-4, to=1-5]
	\arrow["{\cdot\ell=c}", from=1-4, to=2-4]
	\arrow[from=2-1, to=2-2]
	\arrow[from=2-2, to=2-3]
	\arrow[from=2-3, to=2-4]
	\arrow[from=2-4, to=2-5]
    \end{tikzcd}\]
By the Snake Lemma, there is an exact sequence $\ker a \rightarrow \ker (b_1 \oplus b_2) \rightarrow \ker c$. By the algebraic $g$-theorem, both $a$ and $b_2$ are injective when $0\leq i<\lfloor \frac{d-1}{2}\rfloor$. It follows that
    $$K_i(\Delta) = \ker b_1 \hookrightarrow \ker c = K_i(\partial\Delta).$$
    Therefore,  $\dim K_i(\Delta) \leq \dim K_i(\partial\Delta)$ for $0\leq i < \lfloor \frac{d-1}{2}\rfloor$, as desired.
\end{proof}

The last ingredient in the proof of Theorem \ref{thm:neighborly} is the following vanishing result for $K_{j-1}(\partial\Delta)$ when $\partial\Delta$ is $i$-neighborly. 

\begin{lemma}\label{lem:K_eq_0}
If $\partial\Delta$ is an $i$-neighborly simplicial $(d-2)$-sphere with $g_1(\partial\Delta)>0$,
then $K_{j-1}(\partial\Delta)=0$ for all $1\leq j\leq i$.
\end{lemma}
\begin{proof}
   Recall that $S = \kk[x_1, \ldots, x_{n+1}]$, where $[n]$ is the vertex set of $\Delta$, and that $\theta_1, \ldots, \theta_d, \ell \in S_1$ are generic linear forms, with $\Theta = (\theta_1,\ldots, \theta_d)$. Relabeling the vertices if necessary, we may assume that $[m]\subseteq[n]$ is the vertex set of $\partial\Delta$. Then $I_{\partial\Delta}\subset S$ is generated by $x_{m+1},\ldots, x_{n+1}$ together with the monomials corresponding to the nonfaces of $\partial\Delta$.  

Since $\partial\Delta$ is $i$-neighborly, all its nonfaces have cardinality greater than $i$. Thus, in degrees at most $i$, the ideal $I_{\partial\Delta}+\Theta$ agrees with  
$(x_{m+1},\ldots,x_{n+1})+\Theta$. As $\Theta$ is generated by $d$ generic linear forms, we obtain an isomorphism of truncated algebras: 
\[
    A:=\bigoplus_{j=0}^{i}\bigl(S/(I_{\partial\Delta}+\Theta)\bigr)_j
    \cong \bigoplus_{j=0}^{i}\bigl(\kk[y_1,\ldots,y_{m-d}]\bigr)_j=:A'.
    \]
By assumption, $m-d=g_1(\partial\Delta)>0$. This, together with genericity of $\ell$, implies that the image of $\ell$ in $A'_1$ is nonzero, and multiplication by this image is injective from $A'_{j-1}$ to $A'_j$ for every $1\le j\le i$. The same then holds for multiplication by $\ell$ on $A$, so $K_{j-1}(\partial\Delta)=0$ for all $1\le j\le i$.
\end{proof}

With Lemmas~\ref{lem:K_i} and~\ref{lem:K_eq_0} in hand, we can now prove the main theorem of this section. Recall that $g_j(\partial\Delta)=h_j(\Delta)-h_{d-j}(\Delta)$. The relation between the $f$- and $g$-vectors implies that a simplicial $(d-2)$-sphere $\Gamma$ with $g_1(\Gamma)=s$ is neighborly if and only if $g_j(\Gamma)=\binom{s+j-1}{j}$ for $1\le j\le\lfloor(d-1)/2\rfloor$.

\begin{theorem} \label{thm:neighborly}
Let $d\geq 6$ and let $h=(1,h_1, \ldots, h_{d-1}, 0)$ be a nonnegative integer vector. Then $h$ is the $h$-vector of a simplicial $(d-1)$-ball whose boundary is a neighborly sphere if and only if $h = (1,0,\ldots, 0)$ or the following conditions are satisfied:
\begin{enumerate}
\item There exists a nonnegative integer $s$ such that
\begin{equation*}
(1, h_1-h_{d-1}, h_2-h_{d-2}, \ldots, h_{\lfloor\frac{d-1}{2}\rfloor}-h_{d-\lfloor\frac{d-1}{2}\rfloor}) = \left(1, s, \binom{s+1}{2}, \ldots, \binom{s+\lfloor\frac{d-1}{2}\rfloor-1}{\lfloor\frac{d-1}{2}\rfloor}\right), 
\end{equation*}
\item $(1, h_1-1, h_2-h_1, \ldots, h_{\lfloor d/2 \rfloor-1} - h_{\lfloor d/2 \rfloor-2}, \max\{h_{\lfloor d/2 \rfloor} - h_{\lfloor d/2 \rfloor-1},0\})$ is an $M$-sequence, and
\item $h_{\lceil d/2 \rceil} \geq h_{\lceil d/2 \rceil +1} \geq \cdots \geq h_{d-1}$.
\end{enumerate}
\end{theorem}

\begin{proof}
We first prove sufficiency. If $h=(1,0,\ldots,0)$, the $(d-1)$-simplex has $h$-vector $h$ and a neighborly boundary. Otherwise, the sequence in (1) is an $M$-sequence, and conditions (2), (1), and (3) imply, respectively, the three sufficient conditions of Theorem~\ref{thm:Kolins}. Hence there is a simplicial $(d-1)$-ball $\Delta$ with $h(\Delta)=h$. The $g$-vector of $\partial \Delta$ is then given by the sequence in (1), implying that $\partial \Delta$ is neighborly.

Now suppose that $h$ is the $h$-vector of a simplicial $(d-1)$-ball $\Delta$ with vertex set $[n]$ and neighborly boundary, and set $s=g_1(\partial\Delta)$. Then condition (1) follows from the neighborliness of $\partial\Delta$, while condition (3) follows from the Billera--Lee conditions (see Corollary~\ref{cor:fromBilleraLee}). It remains to prove condition (2). We consider the cases $s=0$ and $s>0$ separately.

If $s=0$, then $\partial\Delta$ is the boundary of a $(d-1)$-simplex $\sigma$. If $\Delta=\sigma$, then $h=(1,0,\ldots,0)$, and we are done. Otherwise, $\sigma$ is not a face of $\Delta$, so gluing $\Delta$ and $\sigma$ along their common boundary gives a simplicial $(d-1)$-sphere $\Gamma:=\Delta\cup\sigma$ with $h(\Gamma)=h+(0,\ldots,0,1)$. Then 
\[g(\Gamma)=\bigl(1,h_1-1,h_2-h_1,\ldots,h_{\lfloor d/2\rfloor}-h_{\lfloor d/2\rfloor-1}\bigr),\]
and by the $g$-theorem this is an $M$-sequence, as desired.

If $s>0$, then $\partial\Delta$ is a neighborly simplicial sphere, and Lemma~\ref{lem:K_eq_0} gives $\dim K_{i-1}(\partial \Delta)=0$ for all $1 \leq i \leq \lfloor\frac{d-1}{2}\rfloor$. It then follows from Lemma~\ref{lem:K_i} that $\dim K_{i-1}(\Delta)=0$ for all $1\leq i \leq \lfloor\frac{d-1}{2}\rfloor$. 

Let $A=S/(I_{\Delta}+(\Theta,\ell))$. The sequence 
\[
\bigl(\dim_\kk A_0,\dim_\kk A_1,\ldots,\dim_\kk A_{\lfloor d/2\rfloor}\bigr)
\]
is an $M$-sequence, being an initial segment of the Hilbert function of the standard graded algebra $A$. Since for $i\geq 1$,
$\dim_\kk A_i=h_i(\Delta)-h_{i-1}(\Delta)+\dim_\kk K_{i-1}(\Delta)$, it follows that
\[
\dim_\kk A_i=h_i-h_{i-1}\quad\text{for }1\le i<\lfloor d/2\rfloor,
\qquad
\dim_\kk A_{\lfloor d/2\rfloor}\ge\max\{h_{\lfloor d/2\rfloor}-h_{\lfloor d/2\rfloor-1},0\}.
\]
Decreasing the last entry of an $M$-sequence to a nonnegative integer preserves the $M$-sequence property. Therefore,
\[
\bigl(1,h_1-1,h_2-h_1,\ldots,h_{\lfloor d/2\rfloor-1}-h_{\lfloor d/2\rfloor-2},\max\{h_{\lfloor d/2\rfloor}-h_{\lfloor d/2\rfloor-1},0\}\bigr)
\]
is an $M$-sequence.
\end{proof}

\section{New necessary conditions; $h$-vectors of simplicial $5$- and $6$-balls}\label{sec:main}

In this section, we establish new necessary conditions for an integer vector to be the $h$-vector of a simplicial ball and give a complete characterization of the $h$-vectors of simplicial $5$- and $6$-balls.

Throughout this section, let $S=\kk[x_1,\dots,x_n]$.
We begin with our main result on necessary conditions. Its proof draws on the results about socles and generic initial ideals reviewed in Section~\ref{sec:gin}, as well as Lemma~\ref{lem:K_i} from the previous section. (Passing to $S=\kk[x_1,\dots,x_n]$ does not change the dimensions of the spaces $K_i$ in Definition~\ref{def:K}.)

\begin{theorem}\label{thm:main}
Let $d \geq 6$ and let $\Delta$ be a simplicial $(d-1)$-ball. If $h_2(\Delta)>0$ then there exists an integer $m \geq 0$ such that
\begin{enumerate}
\item $(1, (h_1-m)-1, h_2-(h_1-m), \max\{h_3-h_2, 0\})$ is an $M$-sequence, and
\item $(1, h_1-h_{d-1}, h_2-h_{d-2}, \ldots, h_{\lfloor (d-1)/2 \rfloor}-h_{d-\lfloor (d-1)/2 \rfloor})-(0,m,0,\dots,0)$ is an $M$-sequence.
\end{enumerate}
\end{theorem}

\begin{proof}
Let $m=\dim_{\kk}K_1(\Delta)$, with $K_1(\Delta)$ as in Definition~\ref{def:K}. We will show that this choice of $m$ satisfies conditions (1) and (2).

We first address condition (1). Let $I=\gin(I_\Delta+\Theta)$, $A=S/I$ and $B=S/(I+(x_n))$. Then by Theorem \ref{thm:hilb-CM} and Proposition \ref{prop:gin-hilbert},
\[
\dim_\kk A_i= \dim_\kk \big(S/(I_\Delta+\Theta)\big)_i=h_i(\Delta)
\]
for all $i$.
By Proposition \ref{prop:gin-multiplication}, $\dim_\kk K_i(\Delta)$ is equal to the dimension of the kernel of the map
\[
\rho_i : A_i \stackrel {\cdot x_n}\longrightarrow  A_{i+1}
\]
for all $i$. Thus, for every $i\geq 1$,
\begin{align*}
    \dim_\kk B_i = \dim_\kk A_i - \dim_\kk A_{i-1} + \dim_\kk \ker \rho_{i-1}
    = \dim_\kk A_i - \dim_\kk A_{i-1} + \dim_\kk K_{i-1}(\Delta).
\end{align*}
In particular, the beginning of the Hilbert function of $B$ is
\[
(1,h_1-1,h_2-h_1+m,h_3-h_2+\epsilon)
\]
for some $\epsilon\geq 0$. (Indeed, $h_2>0$ implies $h_1>0$, which, in turn, implies $K_0(\Delta)=0$.) 

In order to obtain condition (1), it suffices to show that $\dim_\kk\Soc(B)_1\geq m$. Indeed, any $m$-dimensional subspace of $\Soc(B)_1$ generates an ideal concentrated in degree $1$. Quotienting $B$ by this ideal produces a standard graded algebra whose Hilbert function begins with $(1,h_1-1-m,h_2-h_1+m,h_3-h_2+\epsilon)$. Condition (1) then follows because decreasing the last entry to $\max\{h_3-h_2,0\}$ preserves the $M$-sequence property.

We prove that $\dim_\kk \Soc(B)_1 \geq m$.
Since $\dim_\kk A_2=h_2(\Delta)>0$ and $I$ is strongly stable, the monomial $x_n^2$ cannot belong to $I$. Therefore $x_n \not \in \ker \rho_1$.
This in particular implies that $\ker \rho_1$ injects into $B_1$ via the natural projection $A \to B$.
By Lemma \ref{sstablesocle},
\[
\ker \rho_1=\{ f \in A_1 : x_n f=0\}=\Soc(A)_1.
\]
The image of $\Soc(A)$ in $B$ lies in $\Soc(B)$, so $\ker \rho_1$ is contained in $\Soc(B)_1$.
Thus,
\[\dim_\kk \Soc (B)_1 \geq \dim_\kk \ker \rho_1 =\dim_\kk K_1(\Delta)=m,\]
where we again use Proposition \ref{prop:gin-multiplication} for the middle equality. This completes the proof of condition~(1).

We now show that $m$ satisfies condition (2). Let $J=I_{\partial\Delta}+\Theta\subseteq S$. We consider $\gin (J)$. By Proposition~\ref{prop:gin-hilbert}, $S/J$ and $S/\gin(J)$ have the same Hilbert function. Since $\partial\Delta$ is a simplicial $(d-2)$-sphere, the algebraic $g$-theorem implies that, for $0\le i\le\lfloor(d-1)/2\rfloor$,
\[
\dim_\kk(S/J)_i=\dim_\kk(S/\gin(J))_i=g_i(\partial\Delta)=h_i-h_{d-i}.
\]
In particular, the Hilbert function of $S/\gin(J)$ coincides with the vector
\[
(1, h_1-h_{d-1}, h_2-h_{d-2}, \ldots, h_{\lfloor (d-1)/2 \rfloor}-h_{d-\lfloor (d-1)/2 \rfloor}).
\]
As in the proof of condition (1), it now suffices to show that $\dim_\kk\Soc(S/\gin(J))_1\ge m$.

Consider the maps
\[
\alpha \colon (S/J)_1 \xrightarrow{\cdot\ell} (S/J)_{2}
\qquad\text{and}\qquad
\beta \colon (S/\gin(J))_1 \xrightarrow{\cdot x_n} (S/\gin(J))_{2}.
\]
By Proposition \ref{prop:gin-multiplication}, $\dim_\kk\ker\beta=\dim_\kk\ker\alpha=\dim_\kk K_1(\partial\Delta)$. 
Since $\gin(J)$ is strongly stable,  Lemma \ref{sstablesocle} gives $\ker \beta= \Soc(S/\gin(J))_1$.
Hence
\[
\dim_\kk\Soc(S/\gin(J))_1 =\dim_\kk\ker\beta
=\dim_\kk K_1(\partial\Delta) \ge\dim_\kk K_1(\Delta)=m,
\]
as desired, where in the last inequality we used Lemma \ref{lem:K_i}.
\end{proof}

Together with Corollary \ref{cor:fromBilleraLee} 
and Kolins’ sufficient conditions, Theorem~\ref{thm:main} settles Kolins’ conjecture and extends the complete characterization of $h$-vectors from simplicial $5$-balls to simplicial $6$-balls.

\begin{theorem}\label{thm:5-and-6-balls}
    Let $d\in\{6,7\}$, and let $h=(1,h_1,\ldots,h_{d-1},0)$ be a nonnegative integer vector. Then $h$ is the $h$-vector of a simplicial $(d-1)$-ball if and only if either $h=(1,h_1,0,\dots,0)$ or there exists an integer $m\geq 0$ such that $h-(0,m,0,\dots,0)$ satisfies Kolins' conditions. 
\end{theorem}

\begin{proof}
    For sufficiency, if $h=(1,h_1,0,\ldots,0)$, then $h$ is the $h$-vector of a stacked ball. Otherwise, $h$ is the $h$-vector of a simplicial ball by Corollary~\ref{cor:Kolins}.

For necessity, let $\Delta$ be a simplicial $(d-1)$-ball with $d\in\{6,7\}$.
If $h_2(\Delta)=0$, then $h(\Delta)=(1,h_1(\Delta),0,\ldots,0)$.
Otherwise, Theorem \ref{thm:main} implies the existence of an integer $m\geq 0$ such that
$h'=h-(0,m,0,\ldots,0)$ satisfies Kolins' first two conditions.
By Corollary \ref{cor:fromBilleraLee}, $h$ satisfies Kolins' third condition;
since $h$ and $h'$ differ only in their second entry, so does $h'$.
\end{proof}

\begin{remark}
Although we have stated the main results of this paper for simplicial balls, they also hold for $\Z/2\Z$-homology balls. The proofs of the necessary conditions apply in this setting, while the sufficiency constructions already produce simplicial balls.
\end{remark}

For simplicial balls of dimension at least $7$, we do not have a plausible conjectural characterization of their $h$-vectors. Motivated by our necessary conditions, we ask the following question.
\begin{question}
Let $\Delta$ be a simplicial $(d-1)$-ball with $h$-vector $h$, and set $s=\lfloor d/2\rfloor-2$. Do there exist nonnegative integers $m_1,\ldots,m_s$ such that the vector $h-(0,m_1,\ldots,m_s,0,\ldots,0)$ is either $(1,0,\ldots,0)$ or satisfies Kolins' conditions?
\end{question}

After the initial version of this paper was posted on arXiv, ChatGPT
(OpenAI, GPT-6.2 Astra) produced an argument that appears to give an
affirmative answer to Question 4.4 in all dimensions, extending the
algebraic method used in this paper.

\section*{AI use statement}
During the early stages of this project, we discussed some mathematical ideas with ChatGPT and Claude, but the final proofs do not rely on those discussions. We also used ChatGPT and Claude for improving grammar and proofreading.

\section*{Acknowledgments}
The first author is partially supported by a graduate fellowship from NSF grant DMS-2246399.
The second author is partially supported by JSPS KAKENHI Grant  JP25K06943. The third author is partially supported by NSF grants DMS-2246399 and DMS-2553961.

\printbibliography

\end{document}